\documentclass[11pt]{amsart}
\newif\ifpdf
    \ifx\pdfoutput\undefined
    \pdffalse 
    \else
    \pdfoutput=1 
    \pdftrue
    \fi

    \ifpdf
    \usepackage[pdftex]{graphicx}
    \else
    \usepackage{graphicx}
    \fi
\usepackage{amsmath,amssymb,amsthm,epsfig}

\newcommand\Z{\mathbb Z}

\newtheorem{theorem}{Theorem}[section]
\newtheorem{proposition}[theorem]{Proposition}
\newtheorem{lemma}[theorem]{Lemma}

\theoremstyle{definition}
\newtheorem{definition}[theorem]{Definition}
\theoremstyle{example}
\newtheorem{example}[theorem]{Example}

\input{epsf.tex}

\begin{document}
\ifpdf
    \DeclareGraphicsExtensions{.pdf, .jpg, .tif}
    \else
    \DeclareGraphicsExtensions{.eps, .jpg}
    \fi

\title{Dead end words in lamplighter groups and other wreath products}
\author{Sean Cleary and Jennifer Taback}
\thanks{The first author acknowledges support from PSC-CUNY grant \#64459-0033}
\thanks{The second author acknowledges support from an NSF-AWM Mentoring Travel Grant and
 NSF grant \#0305441, and would like to thank Bowdoin College for their hospitality during the writing of this paper.}

\begin{abstract}
We explore the geometry of the Cayley graphs of the lamplighter
groups and a
 wide range of wreath products.  We show that
these groups have dead end elements of arbitrary depth with
respect to their natural generating sets.  An element $w$ in a
group $G$ with finite generating set $X$ is a dead end element if
no geodesic ray from the identity to $w$ in the Cayley graph
$\Gamma(G,X)$ can be extended past $w$.  Additionally, we describe
some nonconvex behavior of paths between elements in these Cayley
graphs and seesaw words, which are potential obstructions to these
graphs satisfying the $k$-fellow traveller property.
 \end{abstract}

\maketitle

\section{Introduction}
\label{sec:intro}

The main goal of this paper is to investigate the metric
properties of geodesics and balls in the Cayley graphs of
lamplighter groups $L_n$, with respect to the presentation
$$
L_n = \langle a,t | [t^{i}a t^{-i},t^{j}at^{-j}], a^n
\rangle$$ as well  more general classes of wreath products.
 The lamplighter groups $L_n$ have an easy geometric interpretation,
in which we can view an element of $L_n$ as a set of instructions
for changing the state of a bi-infinite string of light bulbs,
each with $n$ states.  For example, in $L_2$, each light bulb is
merely on or off.

Not much is known about the geometry of the Cayley graph of these
groups in the above presentation, and we begin by describing
normal forms for elements of $L_n$ and discussing a special family
of geodesics in this Cayley graph.  We then describe an
interesting phenomenon which occurs in this Cayley graph, namely
the existence of {\em dead end} elements.  Let $G$ be a group with
a finite generating set $X$, and denote by $\Gamma = \Gamma(G,X)$
the Cayley graph of $G$ with respect to this finite generating
set. An element $w \in G$ is a dead end element if no geodesic ray
from the identity to $w$ in $\Gamma$ can be extended past $w$.

Dead end behavior occurs in a number of settings.  Finite groups
always have dead end elements, and even well-behaved groups such
as the integers can have dead end elements with respect to
contrived generating sets, as discussed in de la Harpe
\cite{delaharpe}.  Examples of infinite groups with dead end
elements with respect to natural generating sets are rare.
Thompson's group $F$ is an example of a torsion-free infinite
group with dead end elements with respect to its natural finite
generating set,
 as described in \cite{ct2}. In Theorem \ref{thm:deadends} below,
 we prove that the lamplighter groups $L_n$ with the generating set
given above have infinite families of dead end elements.

There is a fundamental difference between the dead end elements in
Thompson's group $F$ and the lamplighter groups $L_n$,
characterized by the notion of {\em depth}. The depth of a dead
end element $w$ with $|w| = n$ is the length of the shortest path
from $w$ to a point in $B(n+1)$.  In $F$ the depth of any dead
element is two.  We show in Theorem \ref{thm:deadends} that the
lamplighter groups contain dead end words which require paths of
arbitrary length to leave their containing balls.

We explore other families of words with interesting properties
regarding which generators may decrease their word length.  
We
define a family of {\em seesaw} words with the property that
exactly two generators reduce word length, and furthermore
that that initial choice of length-reducing generator determines a unique
subsequent length-reducing generator for many iterations.

We continue our discussion of the Cayley graph of $L_n$ in Section
\ref{sec:ac} by considering convexity properties. Cannon defines
the notion of {\em almost convexity} in \cite{cannon}.  This
property guarantees algorithmic construction of the ball $B(n+1)$
from the ball $B(n)$ by making it sufficient to consider only a
finite set of possible ways that an element  in $B(n+1)$ can be
obtained from different elements of $B(n)$.  Geometrically, a
group is almost convex with respect to a finite generating set if,
for large enough $n$, any points at distance $2$ in $B(n)$ can be
connected by a path of uniformly bounded length which lies
entirely inside of $B(n)$.

The groups $L_n$ are not finitely presented, and thus it follows
from Cannon \cite{cannon} that the Cayley graphs $\Gamma =
\Gamma(L_n,\{a,t\})$ are not almost convex.  Using the geometry of
the groups $L_n$ as well as our description of geodesics in
$\Gamma$, we give an explicit example of an infinite family of
words which illustrates concretely the failure of the almost
convexity property, and extend these examples to a more general
class of wreath products.

Kapovich \cite{ilya} defines a weaker convexity condition for
Cayley graphs, called $K'(2) $ or {\em minimal almost convexity}.
Kapovich \cite{ilya} and Riley \cite{riley} prove several
consequences of minimal almost convexity; for example, if a group
$G$ is minimally almost convex then Kapovich \cite{ilya} showed
that it is necessarily finitely presented.  It follows that the
lamplighter groups are not minimally almost convex.  We  describe
an explicit family of examples which show concretely why the
groups $L_n$ are not minimally almost convex.

Finally, in Section \ref{sec:genwreath} we describe dead end
behavior in more general wreath products. The lamplighter groups
are algebraically the wreath products $\Z_n \wr \Z$.  The proof of
Theorem \ref{thm:deadends} uses certain properties of this wreath
product structure, as well as metric properties of the groups. We
extend Theorem \ref{thm:deadends} to a larger class of wreath
products in Theorem \ref{thm:deadendswreath} and show
 the existence of `seesaw' words in more general wreath products.

\section{Background}
\label{sec:background}

\subsection{Wreath products}
A wreath product is a standard algebraic construction which is a
special case of a semi-direct product. Given groups $G$ and $H$,
we form the wreath product $G \wr H$ by taking the direct product
of copies of $G$, one for each element of $H$.  The generators of
$G$ act on the copy of $G$ indexed by the identity, and the
generators of $H$ act on the indexing elements of $H$, with the
effect of translating the copy of $G$ indexed by the identity to
the appropriate conjugate copy of $G$.  The direct sum
$\displaystyle \bigoplus_{h \in H} G$ is a subgroup of $G \wr H$
given by conjugate copies of $G$, and there is a natural projection
homomorphism from  $G \wr H$ to $H$ which merely deletes occurrences
of the generators of $G$.

One of the simplest infinite wreath products is $Z_2 \wr \Z$, referred to
as the lamplighter group.   This group has been studied for its
remarkable spectral measure \cite{dickspec, spec2, zukspec}, its
strange random walks \cite{rw1, rw2} and its interest as
an unusual amenable group \cite{lindenstrauss}.

Wreath products often provide interesting examples of algebraic or
geometric results.  For example, Erschler \cite{erschler} proved
that $\Z \wr \Z$ is quasi-isometric to $\Z \wr (\Z \times S)$,
where $S$ is any finite group.  This gave the first example of a
finitely generated, although not finitely presented, solvable
group quasi-isometric to a finitely generated group which is not
virtually solvable.

\subsection{Convexity conditions}

We now define the convexity conditions which we use below.

\begin{definition}
A group $G$ is {\em almost convex$(k)$}, or $AC(k)$ with respect
to a finite generating set $X$ if the following condition holds.
There is a number $N(k)$ so that for all positive integers $n >
N_0$, any two elements $y$ and $z$ in the ball $B(n)$ of radius
$n$ with $d_X(y,z) \leq k$ can be connected by a path of length at
most $N(k)$ which lies entirely in $B(n)$.
\end{definition}

Cannon \cite{cannon} showed that the parameter $k$ is unnecessary;
if $(G,X)$ is $AC(2)$, then it is also $AC(k)$ for all integral $k
> 2$.  Thus we may refer to $(G,X)$ simply as {\em almost convex},
ignoring the parameter $k$.
 Furthermore, if the condition holds
for all finite generating sets $X$ of $G$, we say that $G$ is {\em
almost convex}.

\begin{definition}
A group $G$ is {\em minimally almost convex}, or $MAC$ with
respect to a finite generating set $X$ if the following condition
holds.  Given any two elements $x$ and $y$ in the ball $B(n)$ of
radius $n> N_0$ which are distance $2$ apart, there is a path
between them of length at most $2n-1$ which remains inside $B(n)$.
\end{definition}

Any two points $x, \ y \in B(n)$ can be connected by a path of
length at most $2n$ which remains entirely inside $B(n)$ by
connecting each point to the identity and joining those paths, so
minimal almost convexity is the weakest of a family of weak
convexity conditions described by Kapovich \cite{ilya}.

\section{Normal forms and geodesics in $L_n$ and wreath products}
\label{sec:lengths}
\subsection{Lamplighter groups}

To understand normal forms of elements in $L_n$, we begin with the
case of $L_2 = \Z_2 \wr \Z$, and the standard geometric
interpretation of elements of $L_2$. Consider a bi-infinite string
of light bulbs, each of which has two states,  0 corresponding to
`off' and 1 corresponding to `on', and a cursor which indicates
the current bulb under consideration. A word in the lamplighter
group $L_2$ is a sequence of movements of the cursor and commands
to turn the bulb at the current location of the cursor on or off.
The identity word is represented by all bulbs being off and the
cursor at the origin. Using the presentation given in Section
\ref{sec:intro} above, the generator $t$ moves the cursor to the
right, $t^{-1}$ moves the cursor to the left, and $a=a^{-1}$
changes the state of the bulb at the cursor's current location.
 Thus, a word $w \in L_2$ is represented by a configuration of
bulbs and the final location of the cursor.  This  final cursor location is
easily computed; if $w$ is written in terms of the generators $a$
and $t$ in the presentation above, the final position of the
cursor is just the exponent sum of $t$.

\begin{figure}\includegraphics[width=3in]{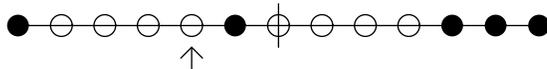}\\
\caption{A typical element of the lamplighter group $L_2$, $w=a_4
a_5 a_6 a_{-1} a_{-6} t^{-2}$. We use a solid circle to represent
a bulb which has been turned on, and the vertical bar denotes the
identity in $\Z$. \label{lampex1}}
\end{figure}

We begin our description of a normal form for elements of $L_2$ by
defining $a_k=t^k a t^{-k}$ for conjugates of $a$.  Beginning with
the bi-infinite string of light bulbs which are all turned off, we
see that $a_n$ moves the cursor to the $n$-th bulb, turns it on,
and returns the cursor to the origin.  It is clear that the $a_n$
all commute. Repeated occurrences of $a_n$ cancel in pairs,
corresponding to turning on a bulb, and later turning it off.

In a product of the generators $a$ and $t$, we can easily move all
instances of $t$ to the end of the word, changing the occurrences
of $a$ to the appropriate $a_k$ to get a word of the form
$$a_{i_1} a_{i_2} \ldots a_{i_k} t^{m}$$
where $i_k \in \Z$.  Since the $a_{i_n}$ commute, we can order
them as we choose. We consider two possible normal forms for a
word $w \in L_2$, separating the word into parts which correspond
to bulbs indexed by negative and non-negative integers:
$$rf(w)=a_{i_1} a_{i_2} \ldots a_{i_k} a_{-j_1} a_{-j_2} \ldots a_{-j_l} t^{m}$$
or
$$lf(w)=a_{-j_1} a_{-j_2} \ldots a_{j_l} a_{i_1} a_{i_2} \ldots a_{i_k}   t^{m}
$$
with $ i_k > \ldots i_2 > i_1 \geq 0 $ and $ j_l > \ldots j_2 >
j_1 > 0 $.

In the `right-first' form, $rf(w)$, the cursor moves first to the
right from the origin, lighting the appropriate bulbs as it moves
toward the bulb with the greatest positive index.  Then the cursor
moves back to the origin and works leftwards, again lighting the
appropriate bulbs in that direction.  Finally, the cursor moves to
its ending location from the leftmost illuminated lamp.

The `left-first' form is similar, but instead of initially moving
to the right, the cursor begins by moving toward the left.  We
mention two specific cases where there are no bulbs illuminated on
one side of the origin, calling these `right-only' or `left-only'
words.

One or possibly both of these normal forms will lead to
minimal-length representation for $w$, depending upon the final
location of the cursor relative to the extreme positive and
negative illuminated bulbs.

We now relate the normal form of an element to its word length in
the group, where word length is  computed with respect to the
generators $\{a,t\}$ from the presentation for $L_2$ given in
Section \ref{sec:intro}.

\begin{definition}
If a word $w \in L_2$ is in either normal form given above, we
define
$$D(w)=k+l+ min\{2 j_l+i_k + | m-i_k|, 2 i_k+j_l+|m+j_l|\}.$$

\end{definition}

We note that geometrically, $D(w)$ is the sum of several
quantities related to the geometric picture of the element $w$:
the number of bulbs which are on, twice the distance of the
furthest bulb from the origin in one direction, the distance of
the furthest bulb from the origin in the other direction and the
distance of the final position of the cursor from the furthest
illuminated bulb in the second direction.

We begin by proving that in $L_2$, the quantity $D(w)$ is exactly
the word length of $w$ with respect to the generating set
$\{a,t\}$.

\begin{proposition}
\label{prop:D} Let $L_2$ be generated by $a$ and $t$, as given
above.  The word length of $w \in L_2$ with respect to the
generating set $\{a,t\}$ is given by $D(w)$.
\end{proposition}

Proposition \ref{prop:D} is proved via the following two lemmas.

\begin{lemma}
\label{lem:upper} The length of a word $w \in L_2$ with respect to
the generating set $\{a,t\}$ is at most $D(w)$.
\end{lemma}

\begin{proof}
We put $w$ into the right-first form listed above and note that by
cancelling adjacent opposite  powers of $t$, we obtain the
expression below for $w$, namely $$w =t ^{-j_1} a t ^{-(j_2-j_1)}
a \ldots t^{-(j_l-j_{l-1})} a t^{j_l+ i_1} a t^{i_2-i_1} a \ldots
t^{i_k-i_{k-1}} a t^{m-i_k}$$ which has $k+l$ occurrences of $a$
and $2 j_l + i_k + |m-i_k|$ occurrences of $t$. Similarly, after
cancellation the left-first form reduces to $$w = t^{ i_1} a
t^{i_2-i_1} a \ldots t^{i_k-i_{k-1}} a t^{-(i_k+j_1} a
 t ^{-(j_2-j_1)} a \ldots  t^{-(j_l-j_{l-1})} a t^{m+j_l}$$ which has $k+l$ occurrences of $a$
 and $2 i_k + j_l + |m+j_l|$
occurrences of $t$.  These bound the length of $w$ above by $D(w)$.
\end{proof}

\begin{lemma}
\label{lem:lower} The length of a word  $w \in L_2$ with respect
to the generating set $\{a,t\}$ is at least $D(w)$.
\end{lemma}

\begin{proof}

To obtain a lower bound on $|w|$, we view $w$ geometrically, as a
collection of light bulbs which are turned on, and a cursor at a
final position $m \in \Z$.  We must relate this picture to the
minimal number of generators $a$ and $t$ needed to create it.  If
$w$ has $n$ light bulbs which are turned on, then a minimal length
representative for $w$ must contain at least $n$ occurrences of
the generator $a=a^{-1}$, each occurrence of which turns on one
bulb.  If $w$ is written in either of the normal forms given
above, then $n = k+l$.

When counting the total occurrences of the generator $t$ in a
minimal representative for $w$, we first consider words with bulbs
illuminated at both positive and negative indices. Recalling that
the exponent sum of $t$ gives the final position of the cursor, we
consider the partial exponent sums on $t$ for a minimal length
representative of $w$.  For instance, at the moment the rightmost
bulb in position $i_k$ is lit, the exponent sum of the generator
$t$ in that prefix must be $i_k$. Similarly, when the leftmost
bulb in position $-j_l$ is lit, the exponent sum of $t$ must be
$-j_l$. Additionally, the total $t$ exponent sum for the entire
word must be $m$.

We consider whether the rightmost or
leftmost illuminated bulb is lit first.  If the rightmost bulb is
illuminated first, then the representative for $w$ has prefixes
with exponent sums of $0,i_k,-j_l$ and $m$ on the $t$ generator.
The total number of $t$ or $t^{-1}$ letters needed to accomplish
this is $i_k+i_k+j_l+|m+j_l|$. Similarly, if the leftmost bulb is
illuminated before the rightmost one, the word must have prefixes
with exponent sums $0,-j_l,i_k$ and then $m$.  This requires at
least $j_l+j_l+i_k+|m-i_k|$ total occurrences of $t$ and $t^{-1}$.
We note that for
 `right-only' and `left-only' words one of $k$ or $l$ will be zero
and the same bounds apply.

Combining the bounds on numbers of $a$ and $t$ generators
appearing in a minimal length representative for $w$, we see that
the lower of these two bounds is exactly $D(w)$.
\end{proof}
It follows immediately that $|w| = D(w)$.

For the lamplighter groups $L_n = \langle a,t | [t^{i}a
t^{-i},t^{j}at^{-j}], a^n \rangle$ with $n>2$, we note that there
are similar normal forms and geodesics. The difference is that the
occurrences of $a$ in the normal form above are replaced by $a^k$,
for $k \in \{ -h, -h+1, \cdots,-1,0,1,2, \cdots h \}$ where $h$ is
the integer part of $\frac{n}{2}$.  When $n$ is even, we omit
$a^{-h}$ to ensure uniqueness, since $a^h=a^{-h}$ in $\Z_{2h}$.

We replace the definition of $D(w)$ given above with the following
more general definition, based on the left and right first normal
forms for elements $w \in L_n$:
$$rf(w)=a^{e_{1}}_{i_1} a^{e_{2}}_{i_2} \ldots a^{e_{k}}_{i_k} a^{f_{1}}_{-j_1} a^{f_{2}}_{-j_2} \ldots a^{f_{l}}_{-j_l} t^{m}$$
or
$$lf(w)=a^{f_{1}}_{-j_1} a^{f_{2}}_{-j_2} \ldots a^{f_{l}}_{-j_l} a^{e_{i}}_{i_1} a^{e_{i}}_{i_2} \ldots a^{e_{k}}_{i_k}   t^{m}
$$
with $ i_k > \ldots i_2 > i_1 \geq 0 $ and $ j_l > \ldots j_2 >
j_1 > 0 $ and $e_i, \ f_j$ lie in  the range $\{ -h,\cdots, h \}$,
as does the exponent $k$ discussed above.  Using these normal
forms we make the following definition.
\begin{definition}
If a word $w \in L_n$ is in normal form as given above, we define
$$D(w)=\Sigma e_{i}+\Sigma f_{j} + min\{2 j_l+i_k + | m-i_k|, 2 i_k+j_l+|m+j_l|\}.$$
\end{definition}
With this definition, $D(w)$ again computes the word length of $w
\in L_n$.

\begin{proposition}
\label{prop:Dn} Let $L_n$ be generated by $a$ and $t$, as given
above.  The word length of $w \in L_n$ with respect to the
generating set $\{a,t\}$ is exactly $D(w)$.
\end{proposition}

The proof is analogous to that of Proposition \ref{prop:D},
additionally taking into account the number of possible states for
each light bulb.

We note that these normal forms by no means give unique minimal
length representatives of elements of $L_n$.  For example, a word
in right-first normal form turns on the bulbs while moving right,
then moves leftwards to the origin back across bulbs already
traversed without changing the state of any bulbs, then continues
moving leftwards, illuminating bulbs while travelling, and then
travels to the final position.  There are potentially many other
minimal-length
 representatives which, for example, move first to the rightmost
 bulb in an `on' state, then turn on bulbs while moving leftwards back to the origin and
so on.  Suppose that a bulb in position $k$ is illuminated, and
the cursor passes position $k$ twice.  Then there is always a
choice whether to illuminate the bulb on the first visit or the
second.  We note that each bulb is visited at most twice, so it is
clear that all minimal length representatives will be of one of
these forms, up to the choices described above.  For each
minimal-length form of an element of $L_2$ (right-first or
left-first) there will be $2^p$ minimal length representatives
where $p$ is the number of `on' bulbs at positions visited twice
during the application of that normal form.

In the case of $L_n$ for $n>2$, there are even more possibilities,
since a bulb which is visited twice may be moved partway to its final state during one
visit to the bulb, and then changed from that state  into its
final state during the subsequent visit.

\subsection{Wreath products}

We note that in wreath products of finitely generated groups with
$\Z$, there is metric structure and normal forms similar to those
described above.

We consider groups of the form $G \wr \Z$, in which we have a
family of minimal length representatives for elements of $G$.  Let
$X$ be a finite set of elements which generate $G$, and let $\Z$
be generated by $t$.  We viewed an element of $L_n$ as a
collection of light bulbs in different positions turned on to a
variety of states; algebraically, this is a collection of elements
of the form $a^{k}$ in $\Z_n$ translated to different conjugate
copies of $\Z_n$ by different powers of the generator $t$ of $\Z$.  In the wreath
products we consider here, we view an element $w \in G \wr \Z$
again as a collection $\{v'_{m}\}$ of elements of $G$, which we
then conjugate into the different copies of $G$ indexed by
elements of $\Z$, together with the final cursor position.
 More precisely, let $v_{m} \in G$ be a minimal
length representative for $v'_{m}$, and $u_{m} =
t^{i_m }v_{m}t^{-i_m}$ the conjugate of $v_{k_m}$ into the
correct conjugate copy of $G$ in the wreath product.  Using this
notation, we again consider right-first and left-first normal
forms for words $w \in G \wr \Z$.  Let
$$rf(w)=u_{i_1} u_{i_2} \ldots u_{i_k} u_{-j_1} u_{-j_2} \ldots u_{-j_l}  t^{m}$$
and
$$lf(w)=u_{-j_1} u_{-j_2} \ldots u_{j_l} u_{i_1} u_{i_2} \ldots u_{i_k}   t^{m}
$$
with $ i_k > \ldots i_2 > i_1 \geq 0 $ and $ j_l > \ldots j_2 >
j_1 > 0 $, and where, as above, $u_m$ is the conjugate by $t^m$ of
a minimal length representative of an element of $G$.  Again we
decide which form is minimal by considering the final position of the
cursor relative to the endpoints.

These normal forms differ from those considered in the
language-theoretic context by Baumslag, Shapiro and Short
\cite{gilbertwreath}; here we choose the order of the $u_k$ to
ensure minimal-length representatives in the wreath product.

We identify the group $G$ with the copy of $G$ in the wreath
product indexed by the identity in $\Z$ and measure the lengths of
elements in any conjugate copy of $G$ by translating them to this
copy.  Thus we have a metric quantity analogous to $D(w)$ defined
above, where we indicate word length in $G$ with respect to the
generating set $X$ by $| \cdot |_G$.

\begin{definition}
Let $w \in G \wr \Z$ be in either normal form given above, and let
$v'_{m}, \ v_{m}, \ i_m$ and $u_{m}$ be as previously defined.
Namely, let $v_{m}$ be a minimal length representative of
$v'_{m} \in G$, and $u_{m}$ a conjugate of $v_{m_i}$ by 
power of $t^{i_m}$.  Then we define
$$D(w)=\sum_{r=1}^{k} {| v_{i_r}|_G+\sum_{s=1}^{l}
| v_{-j_s}|_G}+ min\{2 j_l+i_k + | m-i_k|, 2 i_k+j_l+|m+j_l|\}.$$

\end{definition}

The quantity $D(w)$ is
merely the sum of all the lengths of the elements in the conjugate
copies together with the total number of instances of $t$ needed
to move between those conjugate copies and leave the cursor where
required.

The computation of the word metric in $G \wr \Z$ is analogous to
the case of the lamplighter groups.

\begin{proposition}
\label{prop:genD} Let $G$ be a finitely-generated group. The word
length of $w \in G \wr \Z$ with respect to the generating set $X
\cup \{t\}$ is given by $D(w)$.
\end{proposition}

This proposition is proved by two lemmas identical to Lemmas
\ref{lem:upper} and \ref{lem:lower}, which we omit here.

Again, these right and left normal forms are by no means unique;
there are choices to make analogous to those for the normal forms
in the lamplighter groups.  Let $w \in G \wr \Z$ and consider
those conjugate copies of $G$ containing a nontrivial element
$u_m$ as part of the word $w$.  The choice made in the formulation
of a normal form is between the following three possibilities, as
we imagine the cursor traversing the indexed conjugate copies of
$G$.

\begin{itemize}
\item Represent $u_m$ by a minimal-length form from $G$ entirely
during the first visit to the appropriate conjugate copy of $G$
and do nothing the second visit. \item Do nothing the first visit
and represent $u_m$ by a minimal-length form from $G$ entirely
during the second visit. \item  Represent $u_m$ by a prefix of a
minimal-length form from $G$  during the first visit and the
remaining part of that minimal-length representative during
 the second visit.
\end{itemize}

We note that in more general wreath products, such as $\Z_2 \wr
(\Z \times \Z)$, the computation of the metric becomes much more
complicated because the question of finding minimal length paths
between a collection of points in the integer lattice is more
difficult than in the integer line. In the plane, finding a
minimal length path between a specified set of points is
equivalent to the travelling salesman problem.  Thus, though the
word problem for $\Z_2 \wr (\Z \times \Z)$ is straightforward, the
question of determining minimal-length representatives is NP-hard,
as shown by Parry \cite{parrynp}.

\section{Properties of the Cayley graphs of lamplighter groups}

\subsection{Dead end words in $L_n$}

We begin by formally defining dead end words and their depth.

\begin{definition}
A word $w$ in a finitely generated group $G$ is a {\em dead end
word with respect to a finite generating set $X$} for $G$ if
$|w|=n$ and $|w x| \leq n$ for all generators $x$
 in $X \cup X^{-1}$.
\end{definition}

Such words are called dead end words because a geodesic ray in the
Cayley graph $\Gamma(G,X)$ from the origin to $w$ cannot be
extended beyond $w$. Note that in groups such as $L_2$ where all
relators are of even length, $|w x|$ will be necessarily $n-1$ for
dead end words.

There are different forms of dead end behavior, characterized by {\em depth}.

\begin{definition}
A word $w$ in a finitely generated group $G$ is a {\em dead end
word of depth $k$ with respect to a finite generating set $X$} if
$k$ is the smallest integer with the following property.  If the
word length of $w$ in $n$, then $|w x_1 x_2 \ldots x_l| \leq n$
for $1 \leq l \leq k$ and all choices of generators $x_i \in X
\cup X^{-1}$.
\end{definition}

A dead end word is a place from which is impossible to make
short-term progress away from the identity in the Cayley graph
$\Gamma(G,X)$. The depth of a dead end word reflects how far it
may be necessary to `back up' before being able to move away from
the identity. Dead end words of large depth show that the balls in
the Cayley graph are significantly folded inward upon themselves.
Note that in groups where all relators are of even length, there
are no dead end words of depth one.

We show in \cite{ct2} that there are dead end words of depth two in
Thompson's group $F$ with respect to the standard finite
generating set, but none of depth greater than two.
There is, however, no bound on depth of dead end words
in the lamplighter groups $L_n$.

\begin{theorem}
\label{thm:deadends}
The lamplighter groups $L_n$  contain dead
end words of arbitrary depth  with respect to the generating set
$\{a,t\}$.
\end{theorem}

\begin{figure}\includegraphics[width=3in]{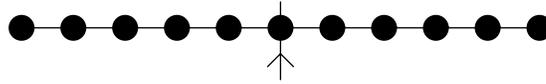}\\
\caption{The word $d_5$, which has the cursor at the origin, and
is a dead end element in $L_2$. A solid circle represent a bulb
which has been turned on. \label{fig:d5}}
\end{figure}

\begin{proof}
We first
prove that $L_n$ contains dead end elements.

We consider a word $d_m$ for $m \geq 1$ in which all of the bulbs
at positions within $m$ of the origin are turned to their setting
furthest from off, and the cursor returns to the origin.  For
example, in $L_2$ the bulbs are merely turned on, and we have the
words $d_m=a_0 a_1 a_2 \ldots a_m a_{-1} a_{-2} \ldots a_{-m} = (a
t) ^m a t^{-m} ( t^{-1} a)^m t^{m}$. See Figure \ref{fig:d5} for
an illustration of $d_5$.  In $L_4$, for example, to put each bulb
into the state farthest from off,
 we apply
the generator $a^2$ to each bulb.  If $h$ is the integer part of $\frac{n}{2}$,
we see that in $L_n$, we want to consider words
$d_m=a_0^h a_1^h a_2^h \ldots a_m^h a_{-1}^h a_{-2}^h \ldots a_{-m}^h =
(a^h t) ^m a^h t^{-m} ( t^{-1} a^h)^m t^{m}$.

It is easily seen that these
words are of length $D_m = 4m+h(2m+1)$ and that the
right-first and left-first normal forms are both of minimal length.

In order for $d_m$ to be a dead end element, each generator of
$L_n$ must decrease the word length of $d_m$.  Writing $d_m$ in
the right-first normal form,  we see that the generator $t^{-1}$
reduces word length. Similarly, writing $d_m$ in the left-first
normal form, we see the generator $t$ must also reduce the word
length of $d_m$. Applying the generator $a$ or $a^{-1}$ turns off
the bulb at the origin in the $n=2$ case.  For even $n>2$, both
$a$ and $a^{-1}$ will decrease the length of the word by making
the bulb at the origin of lower intensity. For odd $n>2$, one of
$a$ and $a^{-1}$ will decrease the length and the other will keep
the length the same.  In all cases except odd $n$, the resulting
words can then be expressed in either normal form using one less
term; in the case of $n$ odd, all of the generators will reduce
length by one except one of either $a$ or $a^{-1}$ which will keep
the length unchanged.

We now show that the elements $d_m$ have depth $m$, providing
examples of dead end words of arbitrary depth in $L_n$.  Let $w$
be a word in which all bulbs at distance greater than $m$ from the
origin are turned off, a collection of bulbs at distance at most
$m$ from the origin are turned on to some possible state, and the
cursor lies somewhere between $-m$ and $m$. It is clear from
$D(w)$ in $L_n$ that the word length of $w$ is at most the word
length of $d_m$.  It follows that an element which is in the ball
$B(k)$ for $k > D_m = |d_m|$ must have a bulb at least distance
$m+1$ from the origin turned on to some state. If $b$ is any such
element, the minimal length of a path from $d_m$ to $b$ is $m+1$,
because in $d_m$, the cursor is at the origin.  Thus we have found
dead end words of depth at least $m$ for arbitrary $m$.
\end{proof}

\subsection{Seesaw words in $L_n$}
\label{seesaw} We again consider the question of when certain
generators decrease the word length of elements of $L_n$.  Dead
end words provide an example of a class of words for which all
generators $\{ a^{\pm 1}, t^{\pm 1}\}$ decrease word length.  We
now describe another class of words with specific conditions on
which generators decrease their word length.  We call these words
{\em seesaw} words, because there is a single generator and its
inverse which reduce word length, and for which successive
applications of that generator also reduce length.

\begin{definition}
A word $w$ in a finitely generated group $G$ with $|w| = n$ is a
{\em seesaw word of swing $k$ with respect to a generator  $g_1$}
in a generating set $X$ if the following conditions hold.
\begin{enumerate}
\item  There is a unique $g \in X$ so that $|wg^{\pm 1}| = |w| -
1$, and for all other $h \in X$, we have $|wh^{\pm 1}| \geq |w|$.

\item Additionally, $|wg^l| = |wg^{l-1}| - 1$ for integral $l \in
[1,k-1]$ and $|wg^{l-1}h^{\pm 1}| \geq |wg^{l-1}|$ for all $h \neq
g \in X$, and the analogous condition for $wg^{-1}$ is also
satisfied.
\end{enumerate}
\end{definition}
Thus $w \in G$ is a seesaw word of swing $k$ if exactly two
generators $g^{\pm 1}$ of $G$ reduce the word length of $w$ and
additionally, $g^{\pm 1}$ is the only generator which reduces the
word length of $wg^{\pm 1}$ for $k-1$ further iterations. Seesaw
words are of interest as they demonstrate one potential difficulty
in choosing geodesic paths which satisfy the `fellow traveller'
property. Namely, the paths $g^{\pm k}$ in the Cayley graph
$\Gamma(G,X)$ from $w$ to $wg^{\pm k}$ are part of minimal length
representatives for $w$ which remain far apart in this graph.

In a finite cyclic group $\Z_{2k}$, the element $k$ is a seesaw
word of swing $k$ with respect to the standard generating set
$\{1\}$.  We prove below that the lamplighter groups contain
seesaw words of arbitrary swing in the generating set used above.
The only other example of seesaw words of this type known to the
authors occur in Thompson's group $F$, and are described in
\cite{ctseesaw}.

\begin{theorem}
\label{thm:seesaw} The lamplighter group  groups $L_n$  contain
seesaw words of arbitrary swing  with respect to the generator
$t$ in the standard generating set $\{a,t\}$.

\end{theorem}

\begin{figure}\includegraphics[width=3in]{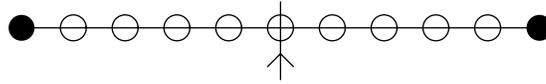}\\
\caption{The seesaw word $w_5 \in L_n$, which has the cursor at
the origin. A solid circle represents a bulb which has been turned
on to any state. \label{lampw5}}
\end{figure}

\begin{proof}
We consider the words $w_n= a^{e_1}_n a^{e_2}_{-n}$, for $e_i \in
\{-h,-h+1, \cdots ,h-1,h\}$ where $h$ is the integer part of
$\frac{n}{2}$.  These words have length $4n+2$. The word $w_5$ is
shown in Figure \ref{lampw5}. For a given $n$, this word consists
of turning on two bulbs to any nontrivial state, a single bulb in
each direction at distance $n$ from the origin and then returning
the cursor to origin.  Right multiplication by both $t$ and
$t^{-1}$ moves the cursor one away from the origin, in each case
shortening the length of $w$. The words $w_n t^k$ are reduced in
length only by further applications of $t$ until the lit bulb is
reached and $a^{-1}$ reduces the length. Similarly, the words $w_n
t^{-k}$ are reduced in length only by further applications of
$t^{-1}$ for $k=1\ldots n-1$, so we have seesaw words of swing
$k$, as desired.

\end{proof}

Since the exponent sum of $t$ in the relators is zero,  words with $t$-exponent sum
differing by $m$ are at least distance $m$ apart.  So these though these two paths from a
seesaw word of swing $k$ are both headed toward the identity, they
immediately separate as fast as possible to distance $2k$.

\section{Convexity properties of some wreath products}
\label{sec:ac}
\subsection{Lamplighter groups are not almost convex}

Cannon \cite{cannon} showed that almost convex groups are finitely
presented, so it is clear that non-trivial wreath products cannot
satisfy the almost convexity condition AC(2) with respect to any
generating set. Kapovich \cite{ilya} showed that groups satisfying
the minimal almost convexity condition are finitely presentable,
so it is clear again that non-trivial wreath products cannot
satisfy the minimal almost convexity condition MAC either.  For
certain wreath products, we can illustrate the failures of these
convexity conditions directly using a natural family of words.

In Example \ref{ex:notac}, we construct a a family of pairs of
words in increasingly large balls in $L_2$ which lie distance two
apart. These words have the property that any path connecting them
which stays inside the ball $B_n$ has a length depending on $n$,
contradicting the definition of almost convexity.

We show directly by example that $L_2$ is not minimally almost
convex in the generating set $\{a,t\}$, and note that the example
extends immediately to $L_n$. We use the relation $a^2 = 1$ to
refer to a bulb as either on or off, but the proof is not affected
by the bulbs having additional ``states" arising from the relation
$a^n = 1$.

We begin by considering the seesaw words $w_n= a_n a_{-n}$, which
were described in Section \ref{seesaw}. These words light bulbs at
distance $n$ from the origin and return the cursor to the origin,
as shown in Figure \ref{lampw5}.
  Right multiplication by both $t$
and $t^{-1}$ moves the cursor one away from the origin, in each
case shortening the length of $w$. The pairs of words we consider
below are always of the form $w_nt$ and $w_nt^{-1}$, each of which
has word length $4n+1$ in the generating set $\{ a,t \}$ of
$L_2$.

\begin{example}
\label{ex:notac} The lamplighter groups $L_n$ do not satisfy
Cannon's almost convexity condition $AC(2)$ with respect to the
generating set $\{a,t\}$.
\end{example}

 We let $\gamma$ be a path in
the Cayley graph $\Gamma = \Gamma(L_2, \{a,t\})$ from $w_nt$ to
$w_nt^{-1}$ which stays entirely inside $B(4n+1)$, where $w_n$ is
defined above. Since the cursor is to the right of the origin in
$w_n$ and to the left of the origin in $w_nt^{-1}$, there must be
at least one prefix of $\gamma$ in which the cursor is at the
origin. Call the first prefix where this occurs $\alpha$.

Let $B_+$ and $B_-$ denote the two light bulbs which are turned on
in both of the words $w_nt$ and $w_nt^{-1}$, with $B_+$ at
position $n$ and $B_-$ at position $-n$.  The following lemma
characterizes the length of words in $L_2$ with these two bulbs
turned on.

\begin{lemma}
\label{lemma:origin}Let $v \in L_2$ be a word with $B_+$ and
$B_-$ turned on, and any additional but possible empty collection
of light bulbs turned on, and the cursor at the origin.  Then $|v|
\geq 4n+2$.
\end{lemma}

\begin{proof}
Given such a word $v$, we can write $v$ in either one of the two
normal forms defined above: the right-first normal form $a_{-j_1}
a_{-j_2} \ldots a_{j_l} a_{i_1} a_{i_2} \ldots a_{i_k}   t^{m}$ or
the left-first normal form $a_{i_1} a_{i_2} \ldots a_{i_k}
a_{-j_1} a_{-j_2} \ldots a_{-j_l}t^{m'}$. In either case, $D(v)
\geq 2+2(i_k+j_l) = 4n + 2 = |w|$, since $i_k \geq n$ and $j_k
\geq n$ and $k$ and $l$ are at least one and $m=m' = 0$.
\end{proof}

To finish this example, we note that since the path $\gamma$ is
contained in $B(4n+1)$, at the point $\alpha$ we must have one of
$B_-$ or $B_+$ turned off, otherwise by Lemma \ref{lemma:origin},
we have $|w \alpha| \geq 4n+2$. The minimum length of a path
turning off one of these bulbs is $n+1$,  depending upon our
original choice for $n$. Thus we see that $L_2$ is not almost convex in the
generating set $\{ a,t \}$.

For lamplighter groups $L_k$ with $k>2$, the situation is
identical, we need only that the bulbs $B_+$ and $B_-$ are in
nontrivial states.

\subsection{Lamplighter groups are not minimally almost convex}
We can use these same examples to show  that $L_n$ is not
minimally almost convex by analyzing the connecting paths more
carefully. Again, since the groups $L_n$ are not finitely
presentable, they are not minimally almost convex by the work of
Kapovich \cite{ilya}; here, we illustrate that fact concretely.
For finitely presented groups, there are not many examples known
of groups which are not minimally almost convex. Belk and Bux
\cite{belkbux} showed that Thompson's group $F$ is not minimally
almost convex, after we showed \cite{ct} that $F$ is not almost
convex.
 Few examples are known of groups which are minimally
almost convex but not almost convex; Elder and Hermiller
\cite{susanmurray} show that the solvable Baumslag-Solitar group
$BS(1,2)$ has this property but some Baumslag-Solitar groups
$BS(1,n)$ with $n>2$ are not minimally almost convex.

\begin{example}
\label{ex:notmac} The lamplighter groups $L_n$ are not minimally
almost convex with respect to the generating set $\{a,t\}$.
\end{example}

Example \ref{ex:notmac} is a continuation of Example
\ref{ex:notac}.  As in Example \ref{ex:notac}, we consider the two
group elements $w_nt$ and $w_nt^{-1}$ which lie in $B(4n+1)$. Let
$\gamma$ be a path from $w_nt$ to $w_nt^{-1}$ which is entirely
contained in $B(4m+1)$. From Lemma \ref{lemma:origin} we know that
at any points along this path which have the cursor at the origin,
at least one of $B_+$ and $B_-$ must be turned off.

At the first point along $\gamma$ with the cursor at the origin,
we must have $B_+$ turned off, otherwise the cursor would be at
the origin with both bulbs illuminated. Since the path goes from
$w_nt$ to $w_nt^{-1}$, we cannot first have $B_-$ turned off, for
that would lead to a point where the cursor is at the origin and
both bulbs are turned on. Thus the length of $\gamma$ is at least
$(n-1)+1=n$, since that is the minimum number of generators
required to move the cursor to position $n$ and turn off the light
bulb.

Since in $w_nt^{-1}$ the bulb $B_+$ is turned on, somewhere along
the path $\gamma$ the bulb $B_+$ must be turned back on. Before
this can happen, $B_-$ must be turned off, otherwise the cursor
would pass through the origin with both bulbs illuminated. Thus
the path $\gamma$ must contain at least an additional $2n+1$
generators, those necessary to move the cursor to the position
$-n$ and turn off the bulb.

We note that $B_+$ must be turned on again before $B_-$ can be
turned on again, since the path $\gamma$ goes from $w_nt$ to
$w_nt^{-1}$, and in $w_nt^{-1}$, the cursor is to the left of the
origin.  Thus $\gamma$ must contain at least an additional $2n+1$
letters which move the cursor to the position $n$ and turn on the
bulb $B_+$, and then an additional $2n+1$ letters which move the
cursor back to position $-n$ and turn on $B_-$.  This brings the
length of $\gamma$ to at least $7n+3$, and we see that to finish
the path there must be at least another $t^{n-1}$, although not
contiguously, returning the cursor to the position $-1$ from the
position $-n$.  Thus the length of $\gamma$ is at least $8n+2$,
which is the sum of the lengths of the original two words,
illustrating that $L_2$ is not minimally almost convex.

\section{More general wreath products}
\label{sec:genwreath}

Some more general wreath products behave similarly to the
lamplighter groups described above.  In the case of dead end
elements, we obtain the following theorem.

\begin{theorem}
\label{thm:deadendswreath} Let $G$ be any finitely generated
non-trivial group which contains dead end elements of depth $d$
with respect to a generating set $X$, and let $\Z$ be generated by
$t$.
 Then $G \wr \Z$ contains dead ends of depth
$d$ with respect to the generating set $X \cup \{t\}$.
 \end{theorem}

\begin{proof}
Let $a$ be a dead end element in $G$ of depth $d$. Define $a_n$ to
be the conjugate $a_n=t^n a t^{-n}$, similar to its definition for
the lamplighter groups above in Theorem \ref{thm:deadends}. The
words $d_n=a_0 a_1 a_2 \ldots a_m a_{-1} a_{-2} \ldots a_{-m}$
will be dead end words of depth $d$ with respect to the generating
set of the wreath product according to an argument identical to
that of Theorem \ref{thm:deadends}.  Namely, right multiplication
by all the generators of $G$ will not increase the word length
since the cursor is at the origin in $d_m$, and $a$ is a dead end
word in $G$; right multiplication by $t$ and $t^{-1}$ will reduce
the length of $d_n$ as well.
\end{proof}

We note that Example \ref{ex:notac} does not rely on the relation
$a^n = 1$ involving the generator $a$ of $L_n$.  It is merely
necessary to have some generator which only affects the second
factor of the wreath product.  Thus, if we can mimic the roles of
the generators $a$ and $t$ in a more general wreath product, the
examples are parallel.

\begin{example}
\label{ex:general} Let $F$ be a finitely generated group
containing an isometrically embedded copy of $\Z$, and $G$ any
finitely generated non-trivial group.  Then $G \wr F$ is neither
almost convex nor minimally almost convex in at least one
generating set.
\end{example}

Again, these groups are not finitely presentable and thus not
almost convex by Cannon \cite{cannon} and not
minimally almost convex by the work of Kapovich \cite{ilya}. We
illustrate this concretely in the following example. Let $a$ be
any generator of $G$, and $t$ the generator of the isometrically
embedded copy of $\Z$ inside of $F$.  Then $\{a,t\}$ can be
completed to a generating set of $G \wr F$.  Again we define $a_n$
to be $t^n a t^{-n}$, and consider the words $w_n=a_n a_{-n}$. The
words $w_n$ give rise to pairs of words  $w_n t^{\pm 1}$ of length
$4n+1$ which are connected by a path of length 2, as in the case
of the lamplighter groups described in Examples \ref{ex:notac} and
\ref{ex:notmac}. The natural analog of Lemma \ref{lemma:origin}
holds; we only use the fact that each word contains a nontrivial
element in the conjugate subgroups $t^n G t^{-n}$ and $t^{-n} G
t^{n}$. Thus we again have concrete examples illustrating the
failure of minimal almost convexity for these groups.

More general wreath products also exhibit the seesaw word
phenomenon described above for the lamplighter groups.

\begin{theorem}
\label{thm:seesawswreath} Let $F$ be a finitely generated group
containing an isometrically embedded copy of $\Z \cong \langle t
\rangle$ with respect to a generating set $Y$, and $G$ any
finitely generated non-trivial group with generating set $X$. Then
$G \wr F$ contains seesaw words of arbitrary swing with respect to
at least one generating set.
 \end{theorem}

\begin{proof}
Let $a \in X$ be any generator of $G$, and $t$ the generator of
the isometrically embedded copy of $\Z$ inside of $F$.  Then
$\{a,t\}$ can be completed to a generating set of $G \wr F$. Again
we define $a_n$ to be $t^n a t^{-n}$, and consider the words of
the form $w_n=a_n a_{-n}$. In these words, the cursor remains
pointed at the conjugate copy of $G$ indexed by the identity.  As
above, right multiplication by any generator of $G$ will increase
the length of $w_n$.  Right multiplication by $t^{\pm 1}$ will
decrease the length of $w_n$ just as it does in the lamplighter
groups considered in Theorem \ref{thm:seesaw}.  Thus these words
are seesaw words in $G \wr F$.  The proof that they have arbitrary
swing is identical to the proof that seesaw words in the
lamplighter groups have this property.
\end{proof}

\bibliographystyle{plain}

\begin{small}
Sean Cleary \\
Department of Mathematics \\
The City College of New York \\
City University of New York \\
New York, NY 10031 \\
E-mail: cleary@sci.ccny.cuny.edu

\medskip

Jennifer Taback\\
Department of Mathematics and Statistics\\
University at Albany\\
Albany, NY 12222\\
E-mail: jtaback@math.albany.edu
\end{small}

\end{document}